\UseRawInputEncoding
\documentclass[14pt]{article}
\usepackage{graphicx}
\usepackage{url}
\usepackage[russian,english]{babel}
\usepackage{amsmath,latexsym,amsthm}
\usepackage{inputenc}
\usepackage{graphics}
\usepackage{graphicx}
\usepackage{amssymb}
\usepackage{color}
\usepackage[colorlinks]{hyperref}
\newtheorem{theorem}{Theorem}[section]
 \newtheorem{definition}[theorem]{Definition}

\newtheorem{proposition}[theorem]{Proposition}
\newtheorem{lemma}[theorem]{Lemma}

\newtheorem{remark}[theorem]{Remark}
\newtheorem{example}[theorem]{Example}
\begin{document}
\newcommand{\tit}[1]{
\begin{center} { {\bf #1 }}  \end{center}}  
\newcommand{\aut}[1]{
\begin{center} {\bf #1 }  \end{center}}           
\tit{ON FUNCTIONAL LIMIT THEOREMS FOR BRANCHING PROCESSES WITH DEPENDENT IMMIGRATION}
\aut{Sadillo Sharipov}

\begin{center}
\emph{V.I.Romanovskiy Institute of Mathematics, Tashkent, Uzbekistan}
\end{center}
\begin{center}
E-mail: sadi.sharipov@yahoo.com
\end{center}

\newcommand{\Abstract}
\noindent \textbf{\Abstract}{ \small ABSTRACT.}
In this paper we consider a triangular array of branching processes with non-stationary immigration. We prove a weak convergence of properly normalized branching processes with immigration to deterministic function under assumption that immigration is rowwise $\psi-$mixing and the offspring mean tends to its critical value 1, immigration mean and variance controlled by regularly varying functions. Moreover, we obtain a fluctuation limit theorem for branching process with immigration when immigration is $m-$dependent where $m$ may tend to infinity with the row index at a certain rate. In this case the limiting process is a time-changed Wiener process. Our results extend and improve the previous known results in the literature.\\

\newcommand{\keywords}
\\

\section{Introduction}
\noindent Let for each $n\in {\mathbb N}$ $\left\{\xi_{k,j}^{\left(n\right)},\, k,j\in {\mathbb N}\right\}$ and $\left\{\varepsilon _{k}^{\left(n\right)},\, k\in {\mathbb N}\right\}$ be two independent families of independent identically distributed (i.i.d) random variables with nonnegative integer values which are defined on a fixed probability space $\left(\Omega,\mathfrak{F}, \mathbb{P} \right)$. The sequence of branching processes with immigration $\left\{X_{k}^{\left(n\right)},\, k\ge 0\right\}$, $n\in {\mathbb N}$ is defined by recursion:
\begin{equation} \label{eq1}
X_{0}^{\left(n\right)}=0, \ \ X_{k}^{\left(n\right)}=\sum_{j=1}^{X_{k-1}^{\left(n\right)}}\xi_{k,j}^{\left(n\right)}+\varepsilon _{k}^{\left(n\right)}, \ \  k,n\in {\mathbb N}.
\end{equation}

Intuitively, for a fixed $n\in {\mathbb N}$, $X_{k}^{\left(n\right)} $ represents the size of $k-$th generation of a
population and $\xi_{k,j}^{\left(n\right)} $ is the number of offspring of the $j-$th individual in the $\left(k-1\right)-$st generation and  $\varepsilon_{k}^{\left(n\right)} $ is the number of immigrants contributing to the $k-$th generation. Note that the assumption $X_{0}^{\left(n\right)}=0$ indicates that the population starts with the offspring of the immigrants in the initial generation.

Assume that for all $n\in {\mathbb N}$, $a_{n}={\mathbb{E}} \xi_{1,1}^{\left(n\right)} <\infty$. For each fixed $n\in {\mathbb N}$, the cases when the offspring mean is less, equal or larger than one are referred to subcritical, critical or supercritical, respectively. The family of  processes \eqref{eq1} is called nearly critical if $a_{n} \to 1$ as $n\to \infty$.

Motivated by results of Rahimov \cite{paper 17}-\cite{paper 18} on a weak convergence of properly normalized scaled process \eqref{eq1} when immigration is independent, it is a natural to ask about generalization of these results for weakly dependent immigration sequences. In this paper, we focus on weak convergence on Skorokhod topology of the properly normalized scaled process \eqref{eq1} when immigration sequence satisfies some mixing conditions.

There have been many research works on limit theorems for \eqref{eq1}. For the first time, Sevastyanov \cite{paper 22} studied continuous-time branching process with immigration when immigration follows a homogenous Poisson process. Then an intensive research has been performed to extend this result to general immigration process. For instance, Nagaev \cite{paper 15} studied the discrete-time process in the case when immigration is a wide sense stationary process and proved that properly normalized branching process with immigration converges in distribution to a gamma distribution. Asadullin and Nagaev \cite{paper 1} considered both discrete and continuous time branching processes with immigration and showed that Nagaev's result \cite{paper 15} still remains true under the more general condition that there exists a random variable $\varepsilon$, such that
${n^{-1}}{\rm E}\left|{\sum\limits_{k=1}^n {\left({{\varepsilon_k}-\varepsilon }\right)}}\right| \to 0$  as $n \to \infty$. For a historical overview on this subject, we refer the reader to the recent survey \cite{paper 21}.

Functional limit theorems (FLT) for \eqref{eq1} have a long history. Many authors considered process \eqref{eq1} and related FLT when immigration sequence is i.i.d. (see \cite{paper 12}, \cite{paper 27}, \cite{paper 26}, \cite{paper 13}-\cite{paper 14}, \cite{paper 5} -\cite{paper 6}, \cite{paper 7}-\cite{paper 8}, \cite{paper 2} and references therein). All these results require independence condition of the immigration process. This constraint can be
relaxed by assuming some reasonable dependence structure: in papers \cite{paper 9}, \cite{paper 11} the immigration sequence is assumed to be a wide--sense stationary process and established a deterministic approximation for \eqref{eq1}. Rahimov\cite{paper 16} considered the case when $\left\{\varepsilon_{k},\, k\in {\mathbb N}\right\}$ is a fixed sequence of independent and non-identically distributed random variables with increasing mean (${\mathbb{E}}\varepsilon_{k}\to \infty$) and proved FLTs for a critical branching process with immigration. Guo and Zhang \cite{paper 4} extended one of the Rahimov result (\cite{paper 16}, Theorem 2.1) to $m-$dependent random variables. Weakly dependent immigration have also recently been studied by Khusanbaev and Sharipov \cite{paper 10} and established a FLT for a critical branching process under a $\phi-$mixing condition. Rahimov \cite{paper 17}-\cite{paper 18} considered process given by \eqref{eq1} and proved diffusion and fluctuation type limit theorems under assumption that the sequence $\left\{\varepsilon_{k}^{\left(n\right)},\, k\in {\mathbb N}\right\}$ is rowwise independent for each $n\geq 1$. In \cite{paper 23}, it is shown that Rahimov's result (\cite{paper 17}, Theorem 1) still holds under assumption that $\left\{\varepsilon_{k}^{\left(n\right)},\, k\in {\mathbb N}\right\}$ satisfies $\phi-$mixing condition (see also \cite{paper 24}-\cite{paper 25}). Further investigations deal with estimation of parameters for process \eqref{eq1} (see \cite{paper 19} -- \cite{paper 20} and references therein).

The purpose of this paper is to extend some FLT in Rahimov \cite{paper 18} and Sharipov \cite{paper 23} to weakly dependent immigration. In Theorem \ref{thm1}, it is proved that scaled random process \eqref{eq1} converges weakly to a deterministic (non-random) function when immigration satisfies $\psi-$mixing condition. Note that Rahimov \cite{paper 18} established the deterministic approximation for \eqref{eq1} a proof of which is based on convergence of random step processes towards a diffusion process. Since in Theorem \ref{thm1} the limiting process is a deterministic, we demonstrate that it can be proven using the classical scheme based on checking the convergence of finite--dimensional distributions and tightness. The next result deals with a FLT for \eqref{eq1}. Under assumption that immigration process is $m-$dependent, we prove a FLT for the fluctuation of \eqref{eq1} (see Theorem \ref{thm2}). To prove Theorem \ref{thm2}, we use the FLT for arrays of martingale differences from \cite{book 3}. Here, $m$ may tend to infinity with the row index at a certain rate. The use of $m-$dependence in Theorem \ref{thm2}, instead of more classical weakly dependence notions relying on the decay rate of mixing coefficients for instance, is motivated by purely technical reason. This reason is that the sequence $\left\{\eta_{k}^{\left(n\right)},k\in {\mathbb N}\right\}$ defined in Lemma \ref{Lemma1} (see Section \ref{section3}) defines $m-$dependent random variables as soon as $\left\{\varepsilon_{k}^{\left(n\right)},\, k\in {\mathbb N}\right\}$ is $m-$dependent for each $n\geq1$.

\indent The rest of the paper is organized as follows. In Section \ref{section2}, we introduce our main assumptions and conditions and provide main results of the paper. Some preliminary results are given in Section \ref{section3}. The proofs of the main theorems are provided in Section \ref{proofs}. For the ease of reading the paper, we provide additional results in the Appendix.

\section{Notations and definitions}\label{section2}

We begin by introducing basic facts.

We recall that a function $f:\left(0,\infty \right)\to \left(0,\infty \right)$ is called regularly varying at infinity if it can be represented in the form $f\left(x\right)=x^{\rho } l\left(x\right)$
where $\rho \in {\mathbb R}$ is called index of regular variation and $l\left(x\right)$ is a slowly varying function.
If a sequence $\left\{f\left(n\right),n\ge 1\right\}$ is regularly varying with exponent $\rho $, we will write $\left\{f\left(n\right),n\ge 1\right\}\in R_{\rho}$.

Assume that for each $n\in {\mathbb N}$, the variables $a_{n}={\mathbb{E}} \xi_{1,1}^{\left(n\right)}$, $b_{n} =Var\left(\xi_{1,1}^{\left(n\right)} \right)$ exist and finite. We also assume that $\alpha \left(n,k\right)={\mathbb{E}} \varepsilon _{k}^{\left(n\right)} <\infty $ and $\beta \left(n,k\right)=Var\left(\varepsilon_{k}^{\left(n\right)} \right)<\infty $ for each $n,k\in {\mathbb N}$. For each $n\in {\mathbb N}$, let $\mathfrak{F}_{k}^{\left(n\right)}$ be $\sigma-$algebra generated by $\left\{X_{0}^{\left(n\right)} ,X_{1}^{\left(n\right)},...,X_{k}^{\left(n\right)} \right\}$.

Denote $A_{n} \left(k\right)={\mathbb{E}} X_{k}^{\left(n\right)}$ and $B_{n}^{2} \left(k\right)=Var\left(X_{k}^{\left(n\right)} \right)$, $1\leq k\leq n$, $n\geq 1$. By \eqref{eq1} and recurrence, we have
$$ A_{n} \left(k\right)=\sum_{j=1}^{k}a_{n}^{k-j} \alpha \left(n,j\right), $$
and from Lemma 1 in \cite{paper 3}, we know that
\begin{equation} \label{eq2}
B_{n}^{2} \left(k\right)=\Delta_{n}^{2} \left(k\right)+\widetilde{\sigma}_{n}^{2} \left(k\right),
\end{equation}
where
$$ \Delta_{n}^{2} \left(k\right)=\frac{b_{n} }{1-a_{n}} \sum_{j=1}^{k}\alpha \left(n,j\right)a_{n}^{k-j-1} \left(1-a_{n}^{k-j} \right), $$
$$ \widetilde{\sigma}_{n}^{2} \left(k\right)=\sigma_{n}^{2} \left(k\right)+2\omega_{n} \left(k\right), $$
$$ \sigma_{n}^{2} \left(k\right)=\sum_{j=1}^{k}\beta \left(n,j\right)a_{n}^{2\left(k-j\right)},\ \  \omega_{n} \left(k\right)=\sum _{j=2}^{k}\sum_{i=1}^{j-1}\operatorname{cov}\left(\varepsilon_{j}^{\left(n\right)},\varepsilon_{i}^{\left(n\right)} \right)a_{n}^{2k-j-i}.$$
Note that in critical case the formula \eqref{eq2} coincides with the formula (2.1) in \cite{paper 4} and when immigration sequence is independent with formula (2.2) in \cite{paper 18}.

In the sequel, we use the following notations: $A\left(n\right)=A_{n} \left(n\right)$, $B^{2} \left(n\right)=B_{n}^{2} \left(n\right)$, $\Delta ^{2} \left(n\right)=\Delta_{n}^{2} \left(n\right)$, $\sigma^{2} \left(n\right)=\sigma_{n}^{2} \left(n\right)$, $\omega \left(n\right)=\omega_{n} \left(n\right)$ in the case when $k=n$. We use the notation $\mathbb{P}-a.s.$ if some relation holds almost surely. $x_{n}\lesssim y_{n}$ means that there exist $C\in \left(0,\infty\right)$, $n_{0}\in \mathbb{N}$ such that $x_{n}\leq C y_{n}$ for all $n\geq n_{0}$, where $\{x_{n},n\geq1\}$ and $\{y_{n},n\geq1\}$ are sequences of real numbers; $\min\left(x,y\right)=x\wedge y$, $x,y\in \mathbb{R}$.

Let $m \geq 1$ be a fixed integer. A sequence $\{\xi_{k}, k \in \mathbb{N}\}$ of random variables is said to be $m-$dependent if the random vectors $\left({{\xi_1},...,{\xi_j}}\right)$ and $\left({{\xi_{j+m+1}},\dots}\right)$ are independent for all $j\geq 1$.
An array $\{\xi_{k}^{\left(n \right)}, k, n \in \mathbb{N}\}$ of random variables is said to be rowwise
$m-$dependent if for every $n \geq 1$, $\{\xi_{k}^{\left(n \right)}, k \in \mathbb{N}\}$ is a sequence of $m-$dependent random
variables. We assume that $m$ depends on the row index $n$ and tends to infinity with appreciate rate (see condition (C5) below).

Now we recall the notion of $\psi-$mixing random variables and arrays of rowwise $\psi-$mixing random variables.

Let $\{\xi_{k}, k \in \mathbb{N}\}$ be a sequence of random variables with zero mean and finite variance and $\mathfrak{F}_a^b$ denotes $\sigma-$algebra generated by the random variables $\xi_i,\,a \leq i \leq b$. Given two $\sigma-$algebras $\mathcal{B}, \mathcal{R}$ in $\mathfrak{F}$, let
$$ \psi \left(\mathcal{B}, \mathcal{R} \right)= \mathop {\sup}\limits_{\scriptstyle A \in \mathcal{B},\scriptstyle B \in \mathcal{R} \hfill \atop
\mathbf{P}\left(A \right)\mathbf{P}\left(B \right)>0 \hfill} \left| \frac{{\mathbf{P}}\left(AB \right)}{{\mathbf{P}}\left(A\right) {\mathbf{P}}\left(B\right)}-1 \right|. $$
Given a sequence $\{\xi_{k}, k \in \mathbb{N}\}$, we associate its sequences of $\psi-$mixing coefficient by letting
$$ \psi\left(n\right)=\mathop {\sup}\limits_{k\geq 1} {\psi}\left(\mathfrak{F}_{1}^{k}, \mathfrak{F}_{k+n}^{\infty} \right). $$
\begin{definition}
A sequence of random variables $\{\xi_{k}, k\in \mathbb{N}\}$  is said to
be $\psi-$mixing if $\psi \left(n \right) \to 0$ as $n \to \infty$.
\end{definition}
A triangular array of random variables $\left\{{{\xi_k^{(n)}},n,k\in \mathbb{N}} \right\}$ is said to be an array of rowwise $\psi-$
mixing random variables if, for every $n\geq 1$, $\left\{{{\xi_k^{(n)}},k\in \mathbb{N}} \right\}$ is a $\psi-$mixing sequence of random variables. Let $\psi_{n} \left(\cdot \right)$ be the mixing coefficient of $\left\{{{\xi_k^{(n)}},k\in \mathbb{N}} \right\}$ for any $n\geq 1$.

The following assumptions are needed in the sequel.

\noindent (C1): There are sequences $\left\{\alpha \left(k\right),k\ge 1\right\}\in R_{\alpha }$ and $\left\{\beta \left(k\right),k\ge 1\right\}\in R_{\beta}$ with $\alpha,\beta \ge 0$, such that, as $n\to \infty$ for each $s\in {\mathbb R}_{+} $,
\begin{equation} \label{eq3}
\mathop{\max }\limits_{1\le k\le ns} \left|\alpha \left(n,k\right)-\alpha \left(k\right)\right|=o\left(\alpha \left(n\right)\right), \ \ \mathop{\max }\limits_{1\le k\le ns} \left|\beta \left(n,k\right)-\beta \left(k\right)\right|=o\left(\beta \left(n\right)\right);
\end{equation}
(C2): $a_{n}=1+an^{-1} +o\left(n^{-1} \right)$ as $n\to \infty $ for some $a\in {\mathbb R}$;

\noindent (C3): $b_{n}=o\left(\alpha \left(n\right)\right)$, $n\to \infty $;

\noindent (C4):  $\alpha \left(n\right)\to \infty $, $\beta \left(n\right)=o\left(n\alpha^{2} \left(n\right)\right)$, $n\to \infty$.

\noindent (C5):  $\alpha \left(n\right)\to \infty $ and $m\to \infty$ such that $m\beta \left(n\right)=o\left(n\alpha \left(n\right)b_{n}\right)$, $n\to \infty$; $ \liminf_{n \to \infty}b_n>0 $.

From the discussion in Rahimov \cite{paper 18}, we know that condition \eqref{eq3} holds if there are sequences of real numbers $x_{n} \to 0$ and $y_{n} \to 0$, $n\to \infty $, such that $\alpha \left(n,k\right)=\alpha \left(k\right)\left(1+x_{n}\right)$ and $\beta \left(n,k\right)=\beta \left(k\right)\left(1+y_{n} \right)$.

In this paper, the symbols $\mathop{\to }\limits^{D}$ and $\mathop{\to}\limits^{P}$ denote the convergence of random functions in Skorokhod space $D\left({\mathbb R}_{+},{\mathbb R}\right)$ and convergence in probability, respectively.

We use the same time change functions as in Rahimov \cite{paper 18} (see (2.3)):
$$ \mu_{\alpha } \left(t\right)=\int_{0}^{t}u^{\alpha} e^{a\left(t-u\right)}du,$$
$$\nu_{\alpha} \left(t\right)=\int_{0}^{t}u^{\alpha}  e^{a\left(t-u\right)} \left(1-e^{a\left(t-u\right)} \right)du,\ \ \lambda_{\beta}\left( t\right)=\int_{0}^{t}u^{\beta} e^{2a\left(t-u\right)}du, $$
\begin{equation} \label{eq4}
\varphi \left(t\right)=\left\{\begin{array}{l}{t^{2+\alpha } ,\, \, \, \, \, a=0,} \\ {\frac{a}{\nu_{\alpha } \left(1\right)} \int _{0}^{t}\mu_{\alpha}\left(u\right)e^{2a\left(t-u\right)}du,\, \, \, a\ne 0. } \end{array}\right.
\end{equation}
Note that $\mu_{\alpha}\left(t\right)=\frac{t^{\alpha+1}}{{\alpha +1}}$ when $a=0$ and  $\mathop{\lim}\limits_{a\to 0} \frac{\nu_{\alpha } \left(t\right)} {a}=\frac{t^{\alpha+2}}{\left(\alpha +1\right)\left(\alpha +2\right)}$.

For $t\in {\mathbb R}_{+}$, define random processes
\begin{equation} \label{eq5}
{X_n}\left(t \right)=\frac{X_{\left[{nt} \right]}^{\left(n \right)} }{A\left(n\right)}, \ \ Z_{n} \left(t\right)=\frac{X_{[nt]}^{\left(n\right)}-A_{n}\left(\left[nt\right]\right)}{B\left(n\right)},\ \ n\ge 1,
\end{equation}
where $\left[x\right]$ denotes the integer part of nonnegative real number $x$.

We denote by $\left\{Y_{n,i}^j\left(k \right),\,k\geq i\right\}$, $1\leq j\leq \varepsilon_{i}^{\left(n\right)}-$the Galton-Watson branching process generated by $j-$th particle arriving at the moment $i$ in the $n-$th series with $Y_{n,i}^j\left(i \right)=1$. From our assumptions it follows that processes $\left\{Y_{n,i}^j\left(k \right),\,k \geqslant i\right\}$, $i,j \geqslant 1$ are independent; $Y_{n,i}^j\left({k+i}\right),\,k\geq0$ has the same distribution as $Y_{n,1}^1\left(k \right),\,k\geq 1$.

Set
$${f_{n,k}}\left(z \right)={ \mathbb{E}}{e^{izY_{n,1}^1\left(k \right)}}, \ \ {\Psi_n}\left(z,t \right)=\mathbb{\mathbb{E}}{e^{izX_{[nt]}^{\left(n \right)}}}, \ \ z\in \mathbb{R}, \ \ t\in \mathbb{R_{+}}.$$
The form of ${\Psi_n}\left(z,t \right)$ is given by relation (3) in \cite{paper 15} (see page 177).

Our first result reads as follows.
\begin{theorem}\label{thm1}
Let $\left\{\varepsilon_{k}^{\left(n\right)},\, n,k\in {\mathbb N}\right\}$ be an array of rowwise $\psi-$mixing random variables satisfying $\sup_{n\geq 1}\sum_{k=1}^{\infty}{\psi_{n}\left(k\right)}<\infty$. If conditions (C1)-(C4) hold, then
\begin{equation} \label{eq7}
X_{n} \left(t\right)\mathop{\to }\limits^{D} \pi_{\alpha} \left(t\right),\ \  n\to \infty
\end{equation}
in Skorokhod space $D\left({\mathbb R}_{+},{\mathbb R}\right)$, where $\pi_{\alpha} \left(t\right)={\mu_{\alpha} \left(t\right)\mathord{\left/ {\vphantom {\mu_{\alpha} \left(t\right) \mu_{\alpha} \left(1\right)}} \right. \kern-\nulldelimiterspace} \mu _{\alpha } \left(1\right)}$, $t\in {\mathbb R}_{+}$.
\end{theorem}
\begin{remark}
Theorem \ref{thm1} extends the corresponding results of \cite{paper 18}, \cite{paper 23} to the case of arrays of rowwise $\psi-$mixing sequence. More precisely, Rahimov \cite{paper 18} obtained \eqref{eq7} when $\left\{\varepsilon_{k}^{\left(n\right)},\, n,k\in {\mathbb N}\right\}$ is independent with conditions (C1)-(C4). In \cite{paper 23}, it is considered the case when $\left\{\varepsilon_{k}^{\left(n\right)},\, n,k\in {\mathbb N}\right\}$ is a $\phi-$mixing with $\sum_{k=1}^{\infty}{\phi^{{1 \mathord{\left/ {\vphantom {1 2}} \right.
\kern-\nulldelimiterspace} 2}}\left(k\right)}<\infty$ and proved \eqref{eq7}.
\end{remark}
It is natural that the examples and Corollaries 2.1-2.3 from Rahimov \cite{paper 18} related to the maximum and the total progeny of the process remain true in the case of dependent immigration. We here provide one more example of application of Theorem \ref{thm1}.
\begin{example}
Let $\xi_{1,1}^{\left(n\right)}$ are Bernoulli random variables with the probability of success $1-an^{-1}$, where $a>0$. Assume that for each $n\in {\mathbb N}$ immigration process $\left\{\varepsilon_{k}^{\left(n\right)},\, k\in {\mathbb N}\right\}$ is rowwise $\psi-$mixing sequence of Poisson distributed random variables with mean $\alpha \left(n,k\right)=\alpha \left(k\right)\left(1+x_{n} \right)$, where $\alpha \left(k\right)\in R_{\alpha }$, $x_{n} \to 0$, $n\to \infty$ and satisfying $\sup_{n\geq 1}\sum_{k=1}^{\infty}{\psi_{n}\left(k\right)}<\infty$. Then  $\beta \left(n,k\right)=\alpha \left(n,k\right)$ and the condition (C1) is fulfilled with $\alpha \left(n\right)=\beta \left(n\right)$. Moreover, it is easily seen that conditions (C2)-(C4) are also satisfied. Hence, we may apply the statement of Theorem \ref{thm1}.
\end{example}
\begin{theorem}\label{thm2} Let for each $n\in {\mathbb N}$, $\left\{\varepsilon_{k}^{\left(n\right)},\, k\in {\mathbb N}\right\}$ be a sequence of $m-$dependent random variables. Assume conditions (C1)-(C3) and (C5) hold and for any $\varepsilon >0$,
\begin{equation}\label{eq8}
{\mathbb{E}} \left(\left(\xi_{k,j}^{\left(n\right)}-a_{n} \right)^{2} I\left\{\left|\xi_{k,j}^{\left(n\right)}-a_{n} \right|>\varepsilon B\left(n\right)\right\}\right) \to 0,\ \ n\to \infty.
\end{equation}
Then
\begin{equation}\label{eq9}
Z_{n} \left(t\right)\, \mathop{\to }\limits^{D} \, W\left(\varphi \left(t\right)\right),\ \ n\to \infty
\end{equation}
in Skorokhod space $D\left({\mathbb R}_{+},{\mathbb R}\right)$, where $\left\{W\left(t\right),t\in {\mathbb R}_{+} \right\}$ is a standard Brownian motion and $\varphi \left(t\right)$ is defined by \eqref{eq4}.
\end{theorem}
\begin{remark}
In fact, condition \eqref{eq8} is the Lindeberg condition on offspring distributions of process \eqref{eq1}. Condition \eqref{eq8} is satisfied if $B^{-\tau}\left(n\right){\mathbb{E}}|\xi_{k,j}^{\left(n\right)}-a_{n}|^{2+\tau} \to 0$ as $n\to \infty$ for some $\tau>0$ (see Remark 2.2 in \cite{paper 18}).
\end{remark}
\begin{remark}
In the case when $\left\{\varepsilon_{k}^{\left(n\right)},\, k\in {\mathbb N}\right\}$ is a sequence of independent random variables and conditions (C1)-(C3) and $\alpha \left(n\right)\to \infty$, $\beta \left(n\right)=o\left(n\alpha\left(n\right)b_{n}\right)$, $n\to \infty$ are fulfilled, Rahimov \cite{paper 18} obtained \eqref{eq9}. We replace condition $\beta \left(n\right)=o\left(n\alpha \left(n\right)b_{n}\right)$ by $m\beta \left(n\right)=o\left(n\alpha\left(n\right)b_{n}\right)$ which is stronger than the latter one, however, it is a natural in the context of rowwise $m-$dependence. Hence, Theorem \ref{thm2} is an improvement for dependent immigration process.
\end{remark}
\begin{remark}
Theorem \ref{thm1} improves the corresponding result of Guo and Zhang \cite{paper 4} to an arrays of nearly critical branching process with immigration.
\end{remark}
\begin{example}
Let $\xi_{1,1}^{\left(n\right)}$ has the following distribution: $\mathbb{P}\left(\xi_{1,1}^{\left(n\right)}=n\right)=\alpha n^{-2}$, $\mathbb{P}\left(\xi_{1,1}^{\left(n\right)}=d_{n}\right)=1/d_{n}^{-1}$, $\mathbb{P}\left(\xi_{1,1}^{\left(n\right)}=0\right)=1-d_{n}^{-1}-\alpha n^{-2}$, where $\alpha \geq 0$ and $d_{n}, n\geq 1$ is a sequence of real numbers such that $d_{n}=O\left(n \right)$, $n \to \infty$. It is easily seen that $a_{n}=1+\alpha n^{-1}$, $b_{n}\sim d_{n}$, $n \to \infty$. Suppose that $\varepsilon_{k}^{\left(n\right)}$ is rowwise $m-$dependent with $m=o\left(n \right)$ as $n \to \infty$ and distributed as $\mathbb{P}\left(\varepsilon_{k}^{\left(n\right)}=[k\ln^{2}k] \right)=\frac{1}{\ln k}\left(1+x_{n}\right)$, $\mathbb{P}\left(\varepsilon_{k}^{\left(n\right)}=0 \right)=1-\frac{1}{\ln k}\left(1+x_{n}\right)$, $n,k\geq 1$. We note that condition (C1) holds with $\alpha\left(n \right)=n\ln n$ and $\beta\left(n \right)=n^{2}\ln^{3} n$ and conditions (C2)-(C3) and (C5) are also satisfied. Condition \eqref{eq8} holds due to fact that for each fixed $\varepsilon >0$ and sufficiently large $n$, the set $\left\{\left|\xi_{k,j}^{\left(n\right)}-a_{n} \right|>\varepsilon B\left(n\right)\right\}$ is empty. Consequently, all conditions of Theorem \ref{thm1} are fulfilled and it remains to apply the statement of Theorem \ref{thm1}.
\end{example}

\section{Auxiliary  Results}\label{section3}
In this section, we give some lemmas which will be used to prove our main results.

We obtain from \eqref{eq1} that
$$ {\mathbb{E}} \left(X_{k}^{\left(n\right)} \left|\mathfrak{F}_{k-1}^{\left(n\right)} \right. \right)=a_{n} X_{k-1}^{\left(n\right)}+{\mathbb{E}} \left(\varepsilon_{k}^{\left(n\right)} \left|\mathfrak{F}_{k-1}^{\left(n\right)} \right. \right). $$
Note that the sequence $\left\{M_{k}^{\left(n\right)},k\ge 1\right\}$ for each $n\in \mathbb{N}$ defined as
$$ M_{k}^{\left(n\right)}:=X_{k}^{\left(n\right)}-{\mathbb{E}} \left(X_{k}^{\left(n\right)} \left|\mathfrak{F}_{k-1}^{\left(n\right)} \right. \right)=X_{k}^{\left(n\right)}-a_{n} X_{k-1}^{\left(n\right)}-{\mathbb{E}} \left(\varepsilon_{k}^{\left(n\right)} \left|\mathfrak{F} _{k-1}^{\left(n\right)} \right.\right) $$
is a martingale difference sequence with respect to the $\sigma-$algebra $\mathfrak{F}_{k}^{\left(n\right)}$, $k\in {\mathbb Z}_{+}$.

On the other hand,
$$ M_{k}^{\left(n\right)}=T_{k}^{\left(n\right)}+N_{k}^{\left(n\right)}, $$
where
\begin{equation}\label{eq10}
T_{k}^{\left(n\right)} =\sum_{j=1}^{X_{k-1}^{\left(n\right)}}\left(\xi_{k,j}^{\left(n\right)}-a_{n} \right), \ \  N_{k}^{\left(n\right)} =\varepsilon_{k}^{\left(n\right)}-{\mathbb{E}} \left(\varepsilon_{k}^{\left(n\right)} \left|\mathfrak{F}_{k-1}^{\left(n\right)} \right. \right).
\end{equation}
The process $Z_{n} \left(t\right)$ given by \eqref{eq5} can be represented as
\begin{equation}\label{eq11}
a_{n}^{-\left[nt\right]} Z_{n} \left(t\right)=Z_{n}^{\left(1\right)} \left(t\right)+Z_{n}^{\left(2\right)} \left(t\right),
\end{equation}
where
$$Z_{n}^{\left(1\right)} \left(t\right)=\frac{1}{B\left(n\right)} \sum_{k=1}^{\left[nt\right]}a_{n}^{-k} M_{k}^{\left(n\right)},$$
$$
Z_{n}^{\left(2\right)} \left(t\right)=\frac{1}{B\left(n\right)} \sum_{k=1}^{\left[nt\right]}a_{n}^{-k} \left({\mathbb{E}} \left(\varepsilon _{k}^{\left(n\right)} \left|\mathfrak{F}_{k-1}^{\left(n\right)} \right. \right)-\alpha \left(n,k\right)\right).
$$
In the following lemma we will prove that $Z_{n}^{\left(2\right)} \left(t\right)$ is asymptotically negligible in $L^{2}-$sense uniformly for all $t\in \left[0,T\right]$, $T>0$.
\begin{lemma}\label{Lemma1}
Assume for each $n\in {\mathbb N}$, $\left\{\varepsilon_{k}^{\left(n\right)},\, k\in {\mathbb N}\right\}$ be a sequence of $m-$dependent random variables. If conditions (C1)-(C3) and (C5) hold, then for each $T>0$,
$$ \sup_{0\leq t \leq T}\left|Z_{n}^{\left(2\right)} \left(t\right)\right|\mathop{\to }\limits^{L^{2} } 0, \ \   n\to \infty. $$
\end{lemma}
\textbf{Proof of Lemma \ref{Lemma1}.} If we denote $\eta_{k}^{\left(n\right)} ={\mathbb{E}} \left(\varepsilon_{k}^{\left(n\right)} \left|\mathfrak{F}_{k-1}^{\left(n\right)} \right. \right)-\alpha \left(n,k\right)$, $k\ge 1$, then for each $n\geq 1$ the sequence $\left\{\eta_{k}^{\left(n\right)},k\in \mathbb{N}\right\}$ also defines $m-$dependent random sequence. Obviously, ${\mathbb{E}} \left(\eta_{k}^{\left(n\right)} \right)^{2}=Var\left({ \mathbb{E}}\left(\varepsilon_{k}^{\left(n\right)}\left|\mathfrak{F}_{k-1}^{\left(n\right)} \right. \right)\right)$. The application of Cauchy--Schwarz inequality and taking into account the inequality $Var\left({ \mathbb{E}} \left(\varepsilon_{k}^{\left(n\right)} \left|\mathfrak{F}_{k-1}^{\left(n\right)} \right. \right)\right)\le \beta \left(n,k\right)$ we have the following bound
$$ {\mathbb{E}} \left(Z_{n}^{\left(2\right)} \left(t\right)\right)^{2}
\le \frac{1}{B^{2} \left(n\right)} \sum_{k=1}^{\left[nt\right]}a_{n}^{-2k} Var\left({\mathbb{E}} \left(\varepsilon_{k}^{\left(n\right)} \left|\mathfrak{F}_{k-1}^{\left(n\right)} \right. \right)\right)+$$
$$ + \frac{2}{B^{2} \left(n\right)} \sum_{k=1}^{\left[nt\right]-1}\sum _{j=k+1}^{\left(k+m\right)\wedge[nt]}a_{n}^{-k} a_{n}^{-j} \sqrt{{\mathbb{E}} \left(\eta_{k}^{\left(n\right)} \right)^{2}} \sqrt{{\mathbb{ E}} \left(\eta_{j}^{\left(n\right)} \right)^{2}}\le $$
\begin{equation}\label{eq13}
\leq \frac{Cm}{B^{2} \left(n\right)} \sum_{k=1}^{\left[nt\right]}a_{n}^{-2k} \beta \left(n,k\right).
\end{equation}
Note that the term \eqref{eq13} can be rewritten as
$$\frac{Cm}{B^{2} \left(n\right)} \sum_{k=1}^{\left[nt\right]}a_{n}^{-2k} \left(\beta \left(n,k\right)-\beta \left(k\right)\right) +\frac{Cm}{B^{2} \left(n\right)} \sum_{k=1}^{\left[nt\right]}a_{n}^{-2k}  \beta \left(k\right), $$
and due to conditions of Theorem \ref{thm2}, we see that the latter relation can be estimated as
$$\frac{mn\beta \left(n\right)}{B^{2} \left(n\right)} \frac{1}{\beta \left(n\right)} \mathop{\max}\limits_{1\le k\le nt} \left|\left(\beta \left(n,k\right)-\beta \left(k\right)\right)\right|\frac{1}{n} \sum_{k=1}^{\left[nt\right]}a_{n}^{-2k} +\frac{mn\beta \left(n\right)}{B^{2} \left(n\right)} \frac{1}{n} \sum_{k=1}^{\left[nt\right]}a_{n}^{-2k}  \to 0, $$
as $n\to \infty$. This ends the proof of Lemma \ref{Lemma1}.

Further, for each $t\in \left[0,T\right]$, $T>0$, let $\Im_{t}^{\left(n\right)}$ be a $\sigma-$field generated by $\left\{Z_{n}^{\left(2\right)} \left(s\right):\, s\le t\right\}$.
Using similar arguments as in the proof of Lemma \ref{Lemma1}, it follows that
$$ {\mathbb{E}} \left(\left|Z_{n}^{\left(2\right)} \left(t+\theta \right)-Z_{n}^{\left(2\right)} \left(t\right)\right|^{2} \left|\Im _{t}^{\left(n\right)} \right. \right) \le {\mathbb{E}} \left(Q_{n} \left(t,\theta \right)\left|\Im_{t}^{n} \right. \right), \ \ \mathbb{P}-a.s. $$
where $\theta \geq 0$ and
$$ Q_{n} \left(t,\theta \right)=\frac{Cm}{B^{2} \left(n\right)} \sum_{j=\left[nt\right]+1}^{\left[nt+n\theta \right]}\left({\mathbb{E}} \left(\varepsilon_{j}^{\left(n\right)} -\alpha \left(n,j\right)\left|\mathfrak{F}_{j-1}^{\left(n\right)} \right. \right)\right)^{2}. $$
Clearly,
\begin{equation}\label{eq14}
Q_{n} \left(t,\theta \right)\le Q_{n} \left(\theta \right):=\mathop{\sup}\limits_{0\le t\le T} Q_{n} \left(t,\theta \right),\ \ \mathbb{P}-a.s.
\end{equation}
and by Lemma \ref{Lemma1}, it yields
$$ \mathop{\sup}\limits_{n\ge 1} \frac{Cm\sum_{k=1}^{\left[nT+n\theta \right]}a_{n}^{-2k} \beta \left(n,k\right) }{B^{2} \left(n\right)} <\infty. $$
By setting
$$G_{n} \left(\theta \right)=\frac{Cm}{B^{2} \left(n\right)} \sum_{k=1}^{\left[nT+n\theta \right]}\left({\mathbb{E}} \left(\varepsilon_{k}^{\left(n\right)}-\alpha \left(n,k\right)\left|\mathfrak{F}_{k-1}^{\left(n\right)} \right. \right)\right)^{2} $$
we obtain
$${\mathbb{E}} \left(G_{n} \left(\theta \right)\right)\le \mathop{\sup}\limits_{n\geq 1} \frac{Cm\sum_{k=1}^{\left[nT+n\theta \right]}a_{n}^{-2k} \beta \left(n,k\right)}{B^{2} \left(n\right)} <\infty $$
and therefore we get $Q_{n} \left(\theta \right)\le G_{n} \left(\theta \right)<\infty, \mathbb{P}-a.s.$
\begin{lemma}\label{Lemma2} Under assumptions of Theorem \ref{thm2}, it holds

1) For each $\eta >0$, $t\ge 0$, there exists $\delta >0$ such that
$$ \mathop{\sup }\limits_{n\ge 1} {\mathbb{P}} \left(Z_{n}^{\left(2\right)} \left(t\right)>\delta \right)\le \eta. $$

2) $$ \mathop{\lim }\limits_{\theta \to 0} \mathop{\sup}\limits_{n\geq 1} {\mathbb{E}} \left(Q_{n} \left(\theta \right)\right)=0.$$
\end{lemma}

\textbf{Proof of Lemma \ref{Lemma2}.}
1) From Chebyshev's inequality and \eqref{eq13}, we get
\begin{equation}\label{eq013}
\mathop{\sup}\limits_{n\ge 1} {\mathbb{P}} \left(Z_{n}^{\left(2\right)} \left(t\right)>\delta \right)\le \frac{1}{\delta^{2}} \mathop{\sup}\limits_{n\ge 1} \frac{Cm\sum_{k=1}^{\left[nt\right]}a_{n}^{-2k} \beta \left(n,k\right)}{B^{2} \left(n\right)}.
\end{equation}
If we denote the right hand side of \eqref{eq013} by $M$ then by Lemma \ref{Lemma1} it follows $M<\infty$. Now, for each fixed $\eta$ and $t$, it remains to choose $\delta $ such that $\delta\geq \left(M/\eta\right)^{{1 \mathord{\left/{\vphantom {1 2}} \right.
\kern-\nulldelimiterspace} 2}}$. Thus, we have proved claim 1).

2) Since $\mathop{\sup}\limits_{0<\theta <1} {\mathbb{E}} \left(Q_{n} \left(\theta \right)\right)\le {\mathbb{E}} \left(G_{n} \left(1\right)\right)\to 0$, $n\to \infty$, then for each $\varepsilon >0$, there exists $N>0$ such that for all $\theta \in \left(0,1\right)$, $n>N$, ${\mathbb{E}} \left(Q_{n} \left(\theta \right)\right)<\varepsilon $. While for all $n\le N$ observing that ${\mathbb{E}} \left(Q_{n} \left(\theta \right)\right)\to 0$ when $\theta \to 0$, one may choose $\theta \in \left(0,1\right)$ such that for $\theta \le \theta_{0}$ and all $k=1,2,...,N$ it yields ${\mathbb{E}} \left(Q_{k} \left(\theta \right)\right)<\varepsilon$. Thus, we have $\mathop{\sup}\limits_{n} {\mathbb{E}} \left(Q_{n} \left(\theta \right)\right)<\varepsilon$ for $\theta \le \theta_{0}$ which proves claim 2). Lemma \ref{Lemma2} is proved.

\begin{lemma}\label{Lemma3}
The set of probability distributions of the process $\left\{Z_{n}^{\left(2\right)} \left(t\right),t\in \mathbb{R}_{+}\right\}$, $n\in \mathbb{N}$, is relatively compact in $D\left({\mathbb R}_{+},{\mathbb R}\right)$.
\end{lemma}
\textbf{Proof of Lemma \ref{Lemma3}.}
From Lemma \ref{Lemma2} we see that conditions (a) and (b) in Ethier and Kurtz \cite{book 2} (see Theorem 8.6.) are fulfilled. Lemma \ref{Lemma3} is proved.\\

Now we provide some lemmas which are taken from Rahimov \cite{paper 18}.
\begin{lemma}\label{Lemma4} Assume $\left\{x\left(n\right),\, n\in {\mathbb N}\right\}\in R_{\rho}$ and condition (C2) holds. Then for each fixed $T>0$ and for all $\rho \in \left[0,\infty \right)$, $\theta \in {\mathbb R}$,
$$ \sup_{0\leq s \leq T} \left|\frac{1}{nx\left(n\right)} \sum_{k=1}^{\left[ns\right]}a_{n}^{k\theta } x\left(k\right)- \int_{0}^{s}t^{\rho}  e^{t\theta a}dt\right| \to 0,  \ \ n\to \infty. $$
\end{lemma}
\begin{lemma}\label{Lemma5} If conditions (C1) and (C2) hold, then uniformly in $s\in \left[0,T\right]$ for each fixed $T>0$,

\begin{enumerate}
\item  $\mathop{\lim }\limits_{n\to \infty } \frac{A_{n} \left(\left[ns\right]\right)}{n\alpha \left(n\right)} =\mu_{\alpha } \left(s\right)$,\ \   $\mathop{\lim }\limits_{n\to \infty } \frac{\sigma_{n}^{2} \left(\left[ns\right]\right)}{n\beta \left(n\right)} =\mu _{\beta } \left(s\right)$,
\item  $\mathop{\lim }\limits_{n\to \infty } \frac{\Delta _{n}^{2} \left(\left[ns\right]\right)}{n^{2} \alpha \left(n\right)b_{n} } =\left\{\begin{array}{l} {{\nu _{\alpha } \left(s\right)\mathord{\left/ {\vphantom {\nu _{\alpha } \left(s\right) a,\, \, \, \, \, \, \, \, \, \, \, \, \, \, \, \, \, \, \, \, \, \, \, \, \, if\, a\ne 0,\, \, \, \, \, \, \, \, \, \, \, \, \, \, \, \, \, \, \, \, \, \, }} \right. \kern-\nulldelimiterspace} a,\, \, \, \, \, \, \, \, \, \, \, \, \, \, \, \, \, \, \, \, \, \, \, \, \, if\, a\ne 0,\, \, \, \, \, \, \, \, \, \, \, \, \, \, \, \, \, \, \, \, \, \, } } \\ {{s^{\alpha +2} \mathord{\left/ {\vphantom {s^{\alpha +2}  \left(\alpha +1\right)\left(\alpha +2\right),\, \, if\, a=0.\, }} \right. \kern-\nulldelimiterspace} \left(\alpha +1\right)\left(\alpha +2\right),\, \, if\, a=0.\, } } \end{array}\right. $
\end{enumerate}
\end{lemma}
\begin{lemma}\label{Lemma6} If conditions (C1) and (C2) hold, then uniformly in $s\in \left[0,T\right]$ for any $\theta \in {\mathbb R}$, $s\in {\mathbb R}_{+} $,
\begin{enumerate}
\item  $\mathop{\lim}\limits_{n\to \infty} \frac{1}{n^{3} \alpha \left(n\right)b_{n}} \sum_{i=1}^{\left[ns\right]}a_{n}^{\theta i} \Delta _{n}^{2} \left(i\right) =\left\{\begin{array}{l} {\left({1\mathord{\left/ {\vphantom {1 a}} \right. \kern-\nulldelimiterspace} a} \right)\int _{0}^{s}e^{u\theta a} \nu _{\alpha } \left(u\right)du,\, \, \, \, \, \, \, \, \, \, \, \, \, \, \, \, \, if\, a\ne 0,\, \, \, \, \, \, \, \, \, \, \, \, \,  } \\ {{s^{\alpha +3} \mathord{\left/ {\vphantom {s^{\alpha +3}  \left(\alpha +1\right)\left(\alpha +2\right)\left(\alpha +3\right),\, \, \, if\, a=0,\, }} \right. \kern-\nulldelimiterspace} \left(\alpha +1\right)\left(\alpha +2\right)\left(\alpha +3\right),\, \, \, if\, a=0,\, } } \end{array}\right. $

\item  $\mathop{\lim }\limits_{n\to \infty } \frac{1}{n^{2} \beta \left(n\right)} \sum_{i=1}^{\left[ns\right]}a_{n}^{\theta i} \sigma _{n}^{2} \left(i\right) =\int_{0}^{s}e^{u\theta a} \mu_{\beta} \left(u\right)du, $

\item  $\mathop{\lim}\limits_{n\to \infty} \frac{1}{n^{2} \alpha \left(n\right)} \sum_{i=1}^{\left[ns\right]}a_{n}^{\theta i} A_{n} \left(i\right) =\int_{0}^{s}e^{u\theta a} \mu_{\alpha } \left(u\right)du $.
\end{enumerate}
\end{lemma}

\begin{lemma}\label{Lemma7} Let $\xi_{k,i}^{\left(n\right)},k,n\in \mathbb{N}$ be a sequence of random variables defined in \eqref{eq1}
and $T_{k}^{\left(n\right)}$ is defined by \eqref{eq10}. Then\\
$1.\ \ {\mathbb{E}} \left(\left(T_{k}^{\left(n\right)} \right)^{2} \left|\mathfrak{F}_{k-1}^{\left(n\right)} \right.\right)=b_{n}X_{k-1}^{\left(n\right)},$\\
$2.\ \ {\mathbb{E}} \left(\left(\sum_{1\leq i \neq j \leq X_{k-1}^{\left(n\right)}}\left(\xi_{k,i}^{\left(n\right)}-a_{n} \right)\left(\xi_{k,j}^{\left(n\right)}-a_{n} \right)\right)^{2}\left|\mathfrak{F}_{k-1}^{\left(n\right)} \right. \right)=2b_{n}^{2}X_{k-1}^{\left(n\right)}\left(X_{k-1}^{\left(n\right)}-1\right).$
\end{lemma}

\section{Proofs of the main theorems}\label{proofs}
In this section, we provide the proofs of our main results.\\
\textbf{Proof of Theorem \ref{thm1}.}
It is known that in order to prove convergence $X_{n} \left(t\right)\mathop{\to }\limits^{D} \pi_{\alpha} \left(t\right)$ as $n\to \infty$, it suffices to verify it for each finite interval. Therefore, let us fix some $T>0$ and denote by $X_{n} \left(t\right)\mathop{\to }\limits^{D\left(T\right)} \pi_{\alpha} \left(t\right)$ the weak convergence of distribution $X_{n} \left(t\right)$ to distribution of $\pi_{\alpha} \left(t\right)$ on $D\left[0,T\right]$.

We only consider the case $a\neq 0$. Then, it is well-known that
\begin{equation}\label{eq16}
{\mathbb{E}}Y_{n,1}^{1}\left(k \right)={a_{n}^k}, \ \  {\mathbb{E}}{\left(Y_{n,1}^{1}\left(k \right) \right)^2}=\frac{{{a_{n}^{k-1}}\left( {{a_{n}^{k}}-1} \right)}}{{a_{n}-1}}{b_{n}}+{a_{n}^{2k}}.
\end{equation}
Since the assumption $b_{n}<\infty$, $n\geq 1$ is equivalent to ${\mathbb{E}}{\left(Y_{n,1}^{1}\left(k \right)\right)^2} < \infty$, we see that the Taylor series expansion is valid for the characteristic function ${f_{n,k}}\left(z \right)$:
\begin{equation}\label{eq6}
{f_{n,k}}\left(z \right)=1+iz{ \mathbb{E}}Y_{n,1}^1\left(k \right) - \frac{{{z^2}}}{2}{\mathbb{E}}{\left({Y_{n,1}^1\left(k \right)} \right)^2} + \frac{{{z^2}}}{2}{\tau_{n,k}}\left(z \right), \ \ k,n \in {\Bbb N},
\end{equation}
where ${\tau_{n,k}}\left(z \right)$ is the remainder term and  $\left| {{\tau_{n,k}}\left(z \right)} \right| \leq 3{\mathbb{E}}{\left( {Y_{n,1}^1\left(k \right)} \right)^2}$, ${\tau_{n,k}}\left(z \right) \to 0$ as $z \to 0$.
Now using decomposition $\ln x=x-1+O\left( {{{\left({x-1}\right)}^2}} \right)$, $x\to 1$ and \eqref{eq16}-\eqref{eq6}, we obtain
\begin{equation}\label{eq17}
\ln \prod\limits_{k=1}^{\left[{nt} \right]} {f_{n,\left[{nt} \right]-k}^{\varepsilon_k^{\left(n \right)}}\left({\frac{z}{{A\left(n \right)}}} \right)}= izI_n^{\left(1 \right)}\left(t \right)-\frac{{{z^2}}}{2}\left({I_n^{\left(2 \right)}\left(t \right)+I_n^{\left(3 \right)}\left(t \right)+I_n^{\left(4 \right)}\left(t \right)} \right),
\end{equation}
where
$$I_n^{\left(1 \right)}\left(t \right)=\frac{{1}}{{A\left(n \right)}}\sum\limits_{k=1}^{\left[{nt} \right]} {a_n^{\left[{nt} \right]-k}\varepsilon_k^{\left(n \right)}},$$
$$I_n^{\left(2 \right)}\left(t \right)=\frac{{1}}{{A^{2}\left(n \right)}}\sum\limits_{k=1}^{\left[{nt} \right]} {\varepsilon_k^{\left(n \right)}{\mathbb{E}}{\left({Y_{n,1}^1\left(\left[{nt} \right]-k \right)} \right)^2}},$$
$$I_n^{\left(3 \right)}\left(t \right)=\frac{{1}}{{A^{2}\left(n \right)}}\sum\limits_{k=1}^{\left[{nt} \right]} {\varepsilon_k^{\left(n \right)}{\tau_{n,\left[{nt} \right]-k}}\left({\frac{z}{{A\left(n \right)}}} \right)},\ \ I_n^{\left(4 \right)}\left(t \right)={\frac{{{1}}}{{A^{2}\left(n \right)}}\sum\limits_{k=1}^{\left[{nt} \right]} {a_n^{2\left({\left[{nt} \right]-k} \right)}\varepsilon_k^{\left(n \right)}}}.$$

We treat each term of relation \eqref{eq17} separately. Let us start with $I_n^{\left(1 \right)}\left(t \right)$. Observe that
\begin{equation}\label{18}
{\mathbb{E}}I_n^{\left(1 \right)}\left(t \right)=\frac{1}{{A\left(n\right)}}\sum\limits_{k=1}^{\left[{nt} \right]} {a_n^{\left[{nt} \right]-k}\left({{\alpha\left(n,k\right)}-{\alpha\left(k\right)}} \right)}+\frac{1}{{A\left(n\right)}}\sum\limits_{k=1}^{\left[{nt} \right]} {a_n^{\left[{nt} \right]-k}{\alpha\left(k\right)}}.
\end{equation}
From Rahimov \cite{paper 18}, we know that under conditions of Theorem \ref{thm1} the first term in \eqref{18} tends to zero and the second term converges to $\pi_{\alpha}\left(t\right)$ uniformly in $t\in \left[0,T\right]$ for each fixed $T>0$.
Thus, ${\mathbb{E}}I_n^{\left(1 \right)}\left(t \right) \to \pi_{\alpha}\left(t\right)$ uniformly in $t\geq 0$.

Now consider the variance of $I_n^{\left(1\right)}\left(t \right)$. Let us consider a sequence $\left\{\zeta_k^{\left(n\right)},k,n\in \mathbb{N} \right\}$ where $\zeta_k^{\left(n\right)}=a_n^{\left[{nt} \right]-k}\varepsilon_k^{\left(n \right)}$, $k,n\in \mathbb{N}$. We see that both sequences
$\left\{\varepsilon_{k}^{\left(n\right)},\, k\in {\mathbb N}\right\}$ and $\left\{\zeta_k^{\left(n\right)},k,n\in \mathbb{N} \right\}$ generate the same $\sigma$-fields. Thus, $\left\{\zeta_k^{\left(n\right)},k,n\in \mathbb{N} \right\}$ is also $\psi-$mixing with ${ \mathbb{E}}\zeta_k^{\left(n\right)}=a_n^{\left[{nt} \right]-k}\alpha\left(n,k\right)$. Consequently, the application of Lemma \ref{Lemma8} (see Appendix) gives us the bound
$$ \operatorname{Var} \left(I_n^{\left(1 \right)}\left(t \right)\right)= \frac{1}{{A^{2}\left(n\right)}}\sum\limits_{k=1}^{\left[{nt}\right]}{\operatorname{Var}\left(\zeta_k^{\left(n\right)}\right)}+ \frac{2}{{A^{2}\left(n\right)}}\sum_{k=2}^{\left[{nt} \right]}\sum_{j=1}^{k-1} \operatorname{cov}\left(\zeta_{k}^{\left(n\right)},\zeta_{j}^{\left(n\right)} \right) $$
$$\leq \frac{1}{{A^{2}\left(n\right)}}\sum\limits_{k=1}^{\left[{nt}\right]}{\operatorname{Var} \left(\zeta_k^{\left(n\right)}\right)}+ \frac{1}{{A^{2}\left(n\right)}}\sum_{k=2}^{\left[{nt} \right]}\sum_{j=1}^{k-1} \psi_{n}\left(k-j \right) \mathbb{E}\left|\zeta_k^{\left(n\right)}\right|\mathbb{E}\left|\zeta_j^{\left(n\right)}\right|  $$
$$\leq \frac{1}{{A^{2}\left(n\right)}}\left(1+\sup_{n\geq 1}\sum_{k=1}^{\infty}\psi_{n}\left(k \right)\right) \sum_{j=1}^{k-1}  \mathbb{E}\left(\zeta_j^{\left(n\right)}\right)^{2}  $$
\begin{equation}\label{eq200}
\leq \frac{C}{{A^{2}\left(n\right)}} \sum_{j=1}^{\left[nt\right]}{a_n^{2\left(\left[{nt}\right]-j\right)}{\beta \left(n,j\right)}}+\frac{C}{{A^{2}\left(n\right)}} \sum_{j=1}^{\left[nt\right]}{a_n^{2\left(\left[{nt}\right]-j\right)}{\alpha^{2} \left(n,j\right)}}.
\end{equation}
It was shown in \cite{paper 18} (see pages 362-363) that under conditions $n\beta\left(n\right)=o\left(A^{2}\left(n\right)\right)$ and (C1) both terms of \eqref{eq200} converge to zero as $n$ tends to infinity.

Combining together the above bounds, we have uniformly in $t\geq 0$
\begin{equation}\label{eq20}
I_n^{\left(1 \right)}\left(t \right)\mathop  \to \limits^P \pi_{\alpha}\left(t\right), \ \  n \to \infty.	
\end{equation}

Consider $I_n^{\left(2 \right)}\left(t\right)$. Applying Lemmas \ref{Lemma4}-\ref{Lemma5} and using the fact that $A\left(n\right)\sim Cn\alpha \left(n\right)$, $n \to \infty$, we derive that
$${\mathbb{E}}I_n^{\left(2\right)}\left(t \right)
= \frac{\Delta_{n}^{2} \left(\left[nt \right]\right)} {A^{2}\left(n\right)}+\frac{1}{{A^{2}\left(n\right)}} \sum_{j=1}^{\left[nt\right]}{a_n^{2\left(\left[{nt}\right]-j\right)}{\alpha \left(n,j\right)}} \to 0.$$
Since $I_n^{\left(2 \right)}\left(t\right) \geq 0$, $\mathbb{P}-a.s.$, we have
\begin{equation}\label{eq21}
I_n^{\left( 2 \right)}\left(t \right) \mathop  \to \limits^P 0, \ \  n \to \infty.
\end{equation}

Using similar arguments as in the proof of $I_n^{\left(2 \right)}\left(t\right)$ and from conditions of Theorem \ref{thm1} and taking into account the inequality $\left| {{\tau_{n,k}}\left(z \right)} \right| \leq 3{\mathbb{E}}{\left({Y_{n,1}^1\left(k \right)} \right)^2}$, one establishes that
$${\mathbb{E}}\left| {I_n^{\left(3 \right)}}\left(t \right) \right| \leq \sum\limits_{k=1}^{\left[{nt} \right]} {\alpha \left(n,k\right)}\left|{\tau_{n,\left[{nt} \right]-k}}\left({\frac{z}{{A\left(n\right)}}} \right)\right|$$
$$ \lesssim  \frac{3}{A^{2}\left(n\right)}\Delta_{n}^{2}\left([nt]\right)+ \frac{3}{A^{2}\left(n\right)}\sum\limits_{k=1}^{\left[{nt}\right]}{a_n^{2\left(\left[{nt} \right]-k\right)}{\alpha \left(n,k\right)}} \to 0,\ \ n \to \infty.$$
Thus by Chebyshev's inequality we deduce that
\begin{equation}\label{eq22}
I_n^{\left(3 \right)}\left(t \right) \mathop  \to \limits^P 0, \ \  n \to \infty.
\end{equation}
Similarly to the proof \eqref{eq21}, we find that $I_n^{\left(4 \right)}\left(t\right)$ converges to zero in probability.

Finally, from the above and \eqref{eq17}, \eqref{eq20}-\eqref{eq22}, we get as $n \to \infty$
$$\ln \prod\limits_{k=1}^{\left[ {nt} \right]} {f_{n,\left[{nt} \right]-k}^{\varepsilon_k^{\left(n \right)}}} \left({\frac{z}{{A\left(n \right)}}} \right)\mathop  \to \limits^P iz \pi_{\alpha}\left(t\right), \ \ t\in \mathbb{R}_{+}.$$
Consequently, the application of Lebesgue dominated convergence theorem gives us
$${\Psi_n}\left({\frac{z}{{A\left(n \right)}},t} \right) \to {e^{iz \pi_{\alpha}\left(t\right)}}, \ \ t\in \mathbb{R}_{+},\ \  n \to \infty.$$
Hence, uniformly with respect to $t \geq 0$
\begin{equation}\label{eq23}
X_{n}\left(t \right) \mathop  \to \limits^P \pi_{\alpha}\left(t\right), \ \  n \to \infty.
\end{equation}
Due to the fact that  the limiting distribution in \eqref{eq23} is degenerate, the convergence of the finite-dimensional distributions follows from \eqref{eq23}. Hence, the finite--dimensional distributions of random process $\left\{{X_n}\left(t \right),t \in \left[{0,T} \right]\right\}$ converge in probability to finite--dimensional distributions of $\left\{\pi_{\alpha} \left(t\right),t \geq 0\right\}$.

It remains to prove the tightness of $\left\{ {X_n}\left(t \right),t \in \left[{0,T} \right]\right\}$.
It is obvious that
$$ {\mathbb{E}}{\left(X_n\left(t \right)-X_n\left(s \right) \right)^2} \leq {K_n^{\left(1 \right)}\left({t,s} \right)+K_n^{\left( 2 \right)}\left({t,s} \right)}, $$
where
$$K_n^{\left(1 \right)}\left({t,s} \right)=\frac{3}{A^{2}\left(n \right)}\left(B_n^{2}\left([nt]\right)+B_n^{2}\left([ns]\right)\right),$$
$$K_n^{\left(2 \right)}\left({t,s} \right)=\frac{3}{A^{2}\left(n \right)}\left(\left(A_n\left([nt]\right)-A_n\left([ns]\right) \right)^2\right).$$
From conditions of Theorem \ref{thm1}, for sufficiently large $n$ and $0\leq s < t \leq T$, we obtain
$$ K_n^{\left(1 \right)}\left({t,s} \right)+K_n^{\left(2 \right)}\left({t,s} \right) = $$
$$= \frac{Cb_n}{{\alpha\left(n\right)}}\left(\nu_{\alpha} \left(t\right)-\nu_{\alpha} \left(s\right) \right)
+\frac{C\beta\left(n\right)}{n\alpha^{2}\left(n\right)}\left(\lambda_{\beta}\left(t\right) -\lambda_{\beta}\left(s\right) \right) + $$
\begin{equation}\label{eq24}
+ \frac{C\max\left(1;e^{at} \right)}{{{n^2}\alpha \left(n\right)}}\sup_{n\geq 1}\sum\limits_{i=1}^\infty {{\psi_{n}}\left(i \right)} \sum_{k=1}^{\left[nt\right]} \beta\left(n,k\right)+T^{2\alpha} \left(t-s \right)^2+o\left(1\right).
\end{equation}

Consequently, from above and \eqref{eq24}, we can argue that
$$ {\mathbb{E}}{\left(X_n\left(t \right)-X_n\left(s \right) \right)^2} \leq CT^{2\alpha} \left(t-s \right)^2. $$
Hence, from Theorem 13.5 in \cite{book 1}, we conclude that the sequence of random processes $\left\{ {X_n}\left(t \right),t \in \left[{0,T} \right]\right\}$ is tight. This ends the proof of Theorem \ref{thm1}.

In the proof of next theorem we use \eqref{eq11} and divide the proof of Theorem \ref{thm2} into two propositions, which
together will imply our result.

\begin{proposition}\label{Proposition1}
Assume for each $n\in {\mathbb N}$, $\left\{\varepsilon_{k}^{\left(n\right)},\, k\in {\mathbb N}\right\}$ is a sequence of $m-$dependent random variables. If conditions (C1)-(C3) and (C5) hold, then
$$Z_{n}^{\left(2\right)} \left(t\right)\mathop{\to }\limits^{D} W\left(\varphi \left(t\right)\right), \ \ n\to \infty $$
in Skorokhod space $D\left({\mathbb R}_{+},{\mathbb R}\right)$, where $\left\{W\left(t\right),t\in {\mathbb R}_{+} \right\}$ is a standard Brownian motion.\\
\end{proposition}
\textbf{Proof of Proposition \ref{Proposition1}.}
First note that since $M_{k}^{\left(n\right)}$ is a martingale difference then the sequence $\left\{U_{k}^{n},k\ge 1\right\}$ where $U_{k}^{n}=a_{n}^{-k}M_{k}^{\left(n\right)}/B\left(n\right)$ for each $n\ge 1$ defines a sequence of martingale differences with respect to the filtration $\mathfrak{F}_{k}^{\left(n\right)} $, $k\in {\mathbb Z}_{+}$. Hence, we need to show that all conditions of Theorem \ref{thm3} are fulfilled (see Appendix). First, we will prove that \eqref{eq40} is satisfied. Observe that
$$ {\mathbb{E}} \left(\left(M_{k}^{\left(n\right)} \right)^{2} \left|\mathfrak{F}_{k-1}^{\left(n\right)} \right. \right)=b_{n} X_{k-1}^{\left(n\right)}+{\mathbb{E}} \left(\left(\varepsilon_{k}^{\left(n\right)} \right)^{2} \left|\mathfrak{F}_{k-1}^{\left(n\right)} \right. \right)-\left({\mathbb{E}} \left(\varepsilon_{k}^{\left(n\right)} \left|\mathfrak{F}_{k-1}^{\left(n\right)} \right. \right)\right)^{2} $$
which yields
\begin{equation}\label{eq25}
\sum_{k=1}^{\left[nt\right]}{\mathbb{E}} \left(\left(U_{k}^{\left(n\right)} \right)^{2} \left|\mathfrak{F}_{k-1}^{\left(n\right)} \right. \right)=J_{n}^{\left(1\right)} \left(t\right)+J_{n}^{\left(2\right)} \left(t\right),
\end{equation}
where
$$ J_{n}^{\left(1\right)} \left(t\right)=\frac{b_{n}}{B^{2}\left(n\right)} \sum_{k=1}^{\left[nt\right]}a_{n}^{-2k} X_{k-1}^{\left(n\right)} ,$$
$$ J_{n}^{\left(2\right)} \left(t\right) = \frac{1}{B^{2} \left(n\right)} \sum_{k=1}^{\left[nt\right]}a_{n}^{-2k} \left({\mathbb{E}} \left(\left(\varepsilon_{k}^{\left(n\right)} \right)^{2} \left|\mathfrak{F}_{k-1}^{\left(n\right)} \right. \right)-\left({\mathbb{E}} \left(\varepsilon_{k}^{\left(n\right)} \left|\mathfrak{F}_{k-1}^{\left(n\right)} \right. \right)\right)^{2} \right).$$

Consider $J_{n}^{\left(1\right)} \left(t\right)$. Since $B^{2} \left(n\right)\sim \Delta^{2} \left(n\right)$, $n\to \infty$ and by using Lemmas \ref{Lemma5}-\ref{Lemma6}, we get uniformly for each $t\in {\mathbb R}_{+}$,
\begin{equation}\label{eq26}
{\mathbb{E}} J_{n}^{\left(1\right)} \left(t\right)=\frac{b_{n}}{B^{2} \left(n\right)} \sum_{k=1}^{\left[nt\right]}a_{n}^{-2k}  A_{n} \left(k-1\right)\to \varphi^\ast\left(t\right), \ \ n\to \infty,
\end{equation}
where $\varphi^\ast\left(t\right)=\left(a/\nu_{\alpha}\left(1\right)\right)\int_0^t\mu_{\alpha}\left(u \right)e^{-2au}du$ if $a\neq 0$ and $\varphi^\ast\left(t\right)=t^{2+\alpha}$ if $a=0$.

Now consider the variance of $J_{n}^{\left(1\right)} \left(t\right)$. It is easily proved that
\begin{equation}\label{eq27}
Var\left(J_{n}^{\left(1\right)} \left(t\right)\right)=R_{n}^{\left(1\right)} \left(t\right)+R_{n}^{\left(2\right)} \left(t\right),
\end{equation}
where
$$ R_{n}^{\left(1\right)} \left(t\right)=\frac{b_{n}^{2} }{B^{4} \left(n\right)} \sum_{k=1}^{\left[nt\right]}a_{n}^{-4k} B_{n}^{2} \left(k-1\right),$$
$$ R_{n}^{\left(2\right)} \left(t\right)=\frac{2b_{n}^{2}}{B^{4} \left(n\right)} \sum_{i=1}^{\left[nt\right]-2}\sum _{j=i+1}^{\left[nt\right]-1}a_{n}^{-2\left(i+j\right)} \operatorname{cov}\left(X_{i}^{\left(n\right)},X_{j}^{\left(n\right)} \right). $$

We will show that $R_{n}^{\left(1\right)} \left(t\right)\to 0$ as $n \to \infty$. With this aim, we first apply the moment inequality for $m-$dependent random variables and obtain
$$ \widetilde{\sigma}^{2}\left(n\right)=Var\left(\sum_{i=1}^{n}a_{n}^{n-k}\left(\varepsilon_{i}^{\left(n\right)}-\alpha\left(n,i\right) \right) \right)$$
$$\leq m\sum_{i=1}^{n}a_{n}^{2\left(i-k\right)} Var\left(\varepsilon_{i}^{\left(n\right)}\right)=m\sigma^{2}\left(n\right).$$
From Lemma \ref{Lemma5} and condition (C5) it follows that $B^{2} \left(n\right)\sim \Delta^{2} \left(n\right)$, $\widetilde{\sigma}^{2}\left(n\right)=o\left(\Delta^{2} \left(n\right) \right)$, $n\to \infty$.
Thus from above and by Lemma \ref{Lemma6}, one can have that
$$
R_{n}^{\left(1\right)} \left(t\right)=\frac{b_{n}^{2} }{B^{4} \left(n\right)} \sum_{k=1}^{\left[nt\right]}a_{n}^{-4k} B_{n}^{2} \left(k-1\right)+\frac{b_{n}^{2}}{B^{4} \left(n\right)} \sum_{k=1}^{\left[nt\right]}a_{n}^{-4k} \widetilde{\sigma}^{2}_{n} \left(k-1\right) $$
\begin{equation}\label{eq28}
\sim \frac{Cn^{3} \alpha \left(n\right)b_{n}^{3} }{n^{4} \alpha^{2} \left(n\right)b_{n}^{2}}+\frac{Cmn^{2} \beta \left(n\right)b_{n}^{2} }{n^{4} \alpha^{2} \left(n\right)b_{n}^{2}} \to 0,\ \ n\to \infty
\end{equation}
where we also used the fact that $m\beta\left(n\right)=o\left(n^{2}\alpha^{2}\left(n \right) \right)$, $n \to \infty$.

Regarding the term $R_{n}^{\left(2\right)} \left(t\right)$, we use the equality
$$ \operatorname{cov}\left(X_{i}^{\left(n\right)},X_{j}^{\left(n\right)} \right)=a_{n}^{j-i} B_{n}^{2} \left(i\right)+\sum_{k=i+1}^{j}\sum _{l=1}^{i}a_{n}^{j-k-i-l} \operatorname{cov}\left(\varepsilon_{k}^{\left(n\right)},\varepsilon_{l}^{\left(n\right)} \right).$$
Unlike the equality (4.8) in Rahimov \cite{paper 18}, the latter formula contains a new term because of dependence of the immigration,
so that we get
$$ R_{n}^{\left(2\right)} \left(t\right)=\frac{2b_{n}^{2}}{B^{4} \left(n\right)} \sum_{i=1}^{\left[nt\right]-2}\sum _{j=i+1}^{\left[nt\right]-1}a_{n}^{-2\left(i+j\right)} a_{n}^{j-i} B_{n}^{2} \left(i\right)+ $$
\begin{equation}\label{eq29}
+\frac{2b_{n}^{2}}{B^{4} \left(n\right)} \sum_{i=1}^{\left[nt\right]-2}\sum_{j=i+1}^{\left[nt\right]-1}a_{n}^{-2\left(i+j\right)} \sum _{k=i+1}^{j}\sum_{l=1}^{i}a_{n}^{j-k-i-l} \operatorname{cov}\left(\varepsilon_{k}^{\left(n\right)},\varepsilon_{l}^{\left(n\right)}\right)  .
\end{equation}

The first term on the right hand-side of \eqref{eq29} can be estimated as
$$ \frac{2b_{n}^{2}}{B^{4} \left(n\right)} \sum_{i=1}^{\left[nt\right]-2}a_{n}^{-3i} B_{n}^{2} \left(i\right) \sum_{j=1}^{\left[nt\right]-1}a_{n}^{-j}\preceq$$
\begin{equation}\label{eq30}
 \preceq\frac{Cn^{3} \alpha \left(n\right)b_{n}^{3} }{n^{4} \alpha^{2} \left(n\right)b_{n}^{2}}+\frac{Cmn^{3}\beta\left(n\right)b_{n}^{2} }{n^{4} \alpha^{2} \left(n\right)b_{n}^{2}} \to 0, \ \ n\to \infty
\end{equation}
where we used conditions (C3) and (C5).

The second term on the right-hand side of \eqref{eq29} can be bounded as
$$ \frac{2b_{n}^{2}}{B^{4} \left(n\right)} \sum_{i=1}^{\left[nt\right]-2}\sum_{j=i+1}^{\left[nt\right]-1}\sum_{k=i+1}^{j}\sum _{l=1}^{i}a_{n}^{-3i-j-k-l} \operatorname{cov}\left(\varepsilon_{l}^{\left(n\right)},\varepsilon_{k}^{\left(n\right)} \right)  $$
$$ \le \frac{2b_{n}^{2}}{B^{4} \left(n\right)} \sum_{i=1}^{\left[nt\right]-2}\sum_{j=i+1}^{\left[nt\right]-1}a_{n}^{-3i-j} \sum _{k=2}^{\left[nt\right]}\sum_{l=1}^{k-1}a_{n}^{-k-l} \operatorname{cov} \left(\varepsilon_{l}^{\left(n\right)},\varepsilon_{k}^{\left(n\right)} \right) $$
\begin{equation}\label{eq31}
\le \frac{Cmb_{n}^{2}}{B^{4} \left(n\right)} \sum_{i=1}^{\left[nt\right]-2}\sum_{j=i+1}^{\left[nt\right]-1}a_{n}^{-3i-j} \sum _{k=2}^{\left[nt\right]}a_{n}^{-2k} \beta\left(n,k\right)  \to 0, \ \ n\to \infty.
\end{equation}
From \eqref{eq29}-\eqref{eq31} we may deduce $R_{n}^{\left(2\right)} \left(t\right)\to 0$, $n\to \infty $. Thus, by  \eqref{eq27}-\eqref{eq28}, we infer that
\begin{equation}\label{eq32}
Var\left(J_{n}^{\left(1\right)} \left(t\right)\right)\to 0,\ \  n\to \infty.
\end{equation}
Recalling \eqref{eq25}-\eqref{eq26}, we get
\begin{equation}\label{eq33}
J_{n}^{\left(1\right)} \left(t\right)\mathop{\to }\limits^{P} \varphi^\ast \left(t\right),\ \  n\to \infty.
\end{equation}
We now prove that $J_{n}^{\left(2\right)} \left(t\right)\mathop{\to }\limits^{P} 0$, $n\to \infty $. Indeed, note that
$$ {\mathbb{E}} J_{n}^{\left(2\right)} \left(t\right)=\frac{1}{B^{2} \left(n\right)} \sum_{k=1}^{\left[nt\right]}a_{n}^{-2k} \left({\mathbb{E}} \left(\varepsilon_{k}^{\left(n\right)} \right)^{2} -\left(\alpha \left(n,k\right)\right)^{2}\right) $$
$$ +\frac{1}{B^{2} \left(n\right)} \sum_{k=1}^{\left[nt\right]}a_{n}^{-2k} \left(\left(\alpha \left(n,k\right)\right)^{2}-{\mathbb{E}} \left({\mathbb{E}} \left(\varepsilon_{k}^{\left(n\right)} \left|\mathfrak{F}_{k-1}^{\left(n\right)} \right. \right)\right)^{2}\right)=$$
$$
=\frac{1}{B^{2} \left(n\right)} \sum_{k=1}^{\left[nt\right]}a_{n}^{-2k} \left(Var\left(\varepsilon_{k}^{\left(n\right)} \right)-Var\left({\mathbb{E}} \left(\varepsilon_{k}^{\left(n\right)} \left|\mathfrak{F}_{k-1}^{\left(n\right)} \right. \right)\right)\right)
$$
\begin{equation}\label{eq34}
\le \frac{1}{B^{2} \left(n\right)} \sum_{k=1}^{\left[nt\right]}a_{n}^{-2k} \beta \left(n,k\right).
\end{equation}
In order to bound the term \eqref{eq34}, we first rewrite it as
$$ \frac{1}{B^{2} \left(n\right)} \sum_{k=1}^{\left[nt\right]}a_{n}^{-2k} \left(\beta \left(n,k\right)-\beta \left(k\right)\right) +\frac{1}{B^{2} \left(n\right)} \sum_{k=1}^{\left[nt\right]}a_{n}^{-2k} \beta \left(k\right) $$
which due to conditions of Theorem \ref{thm2} and Lemma \ref{Lemma4},
$$ \frac{1}{\beta \left(n\right)} \mathop{\max }\limits_{1\le k\le nt} \left|\beta \left(n,k\right)-\beta \left(k\right)\right|\frac{n\beta \left(n\right)}{B^{2} \left(n\right)} \frac{1}{n} \sum_{k=1}^{\left[nt\right]}a_{n}^{-2k} +\frac{n\beta \left(n\right)}{B^{2} \left(n\right)} \frac{1}{n} \sum_{k=1}^{\left[nt\right]}a_{n}^{-2k} \to 0, $$
which proves
\begin{equation}\label{eq35}
J_{n}^{\left(2\right)} \left(t\right)\mathop{\to }\limits^{P} 0,\ \ n\to \infty.
\end{equation}
Consequently, collecting \eqref{eq33} with \eqref{eq35}, we obtain
\begin{equation} \label{eq035}
\sum_{k=1}^{\left[nt\right]}{\mathbb{E}} \left(\left(U_{k}^{n} \right)^{2} \left|\mathfrak{F}_{k-1}^{n} \right. \right) \mathop{\to }\limits^{P} \varphi^{\ast} \left(t\right), \ \  n\to \infty.
\end{equation}

We turn now to the proof of \eqref{eq39}.
Using a simple inequality $\left(a+b\right)^{2} \le 2\left(a^{2}+b^{2} \right)$, where $a,b\in \mathbb{R}$, we find that
$$ L\left(n\right):=\frac{1}{B^{2} \left(n\right)} \sum_{k=1}^{\left[nt\right]}a_{n}^{-2k}{\mathbb{E}} \left(\left(M_{k}^{\left(n\right)} \right)^{2} I\left\{\left|a_{n}^{-k} M_{k}^{\left(n\right)} \right|>\varepsilon B\left(n\right)\right\}\left|\mathfrak{F}_{k-1}^{\left(n\right)} \right. \right)$$
$$\le 2\left(L_{1} \left(n\right)+L_{2} \left(n\right)\right), $$
where
$$ L_{1} \left(n\right)=\frac{1}{B^{2} \left(n\right)} \sum_{k=1}^{\left[nt\right]}a_{n}^{-2k} {\mathbb{E}} \left(\left(T_{k}^{\left(n\right)} \right)^{2} I\left\{\left|a_{n}^{-k} M_{k}^{\left(n\right)} \right|>\varepsilon B\left(n\right)\right\}\left|\mathfrak{F}_{k-1}^{\left(n\right)} \right. \right),$$
$$ L_{2} \left(n\right)=\frac{1}{B^{2} \left(n\right)} \sum_{k=1}^{\left[nt\right]}a_{n}^{-2k} {\mathbb{E}} \left(\left(N_{k}^{\left(n\right)} \right)^{2} I\left\{\left|a_{n}^{-k} M_{k}^{\left(n\right)} \right|>\varepsilon B\left(n\right)\right\}\left|\mathfrak{F}_{k-1}^{\left(n\right)} \right. \right).$$
First we estimate $L_{1} \left(n\right)$. Note that for any random variables $X$ and $Y$, and for all $\varepsilon >0$, one has
\begin{equation} \label{eq36}
I\left\{\left|X+Y\right|>\varepsilon \right\}\le I\left\{\left|X\right|>{\varepsilon \mathord{\left/ {\vphantom {\varepsilon  2}} \right. \kern-\nulldelimiterspace} 2} \right\}+I\left\{\left|Y\right|>{\varepsilon \mathord{\left/ {\vphantom {\varepsilon  2}} \right. \kern-\nulldelimiterspace} 2} \right\}.
\end{equation}
The application of \eqref{eq36} gives us the bound
$$ L_{1} \left(n\right)\le L_{1,1} \left(n\right)+L_{1,2} \left(n\right),$$
where
$$ L_{1,1} \left(n\right)=\frac{1}{B^{2} \left(n\right)} \sum_{k=1}^{\left[nt\right]}a_{n}^{-2k} {\mathbb{E}} \left(\left(T_{k}^{\left(n\right)} \right)^{2} I\left\{\left|a_{n}^{-k} T_{k}^{\left(n\right)} \right|>{\varepsilon B\left(n\right)\mathord{\left/ {\vphantom {\varepsilon B\left(n\right) 2}} \right. \kern-\nulldelimiterspace} 2} \right\}\left|\mathfrak{F} _{k-1}^{\left(n\right)} \right. \right), $$
$$ L_{1,2} \left(n\right)=\frac{1}{B^{2} \left(n\right)} \sum_{k=1}^{\left[nt\right]}a_{n}^{-2k} {\mathbb{E}} \left(\left(T_{k}^{\left(n\right)} \right)^{2} I\left\{\left|a_{n}^{-k} N_{k}^{\left(n\right)} \right|>{\varepsilon B\left(n\right)\mathord{\left/ {\vphantom {\varepsilon B\left(n\right) 2}} \right. \kern-\nulldelimiterspace} 2} \right\}\left|\mathfrak{F} _{k-1}^{\left(n\right)} \right. \right). $$
From Rahimov \cite{paper 18} (see pages 369-370), it is known that
\begin{equation}\label{eq036}
L_{1,1} \left(n\right)\mathop{\to}\limits^{P} 0,\ \ n \to \infty.
\end{equation}
The proof of \eqref{eq036} is omitted since the proof also remains true in our context.

Further, by using Lemma \ref{Lemma7} and taking into account the independence between the offspring and immigration sequences, we get
$$ L_{1,2} \left(n\right)\le \frac{b_{n} }{B^{2} \left(n\right)} \sum_{k=1}^{\left[nt\right]}a_{n}^{-2k} X_{k-1}^{\left(n\right)}{\mathbb{E}} \left(I\left\{\left|a_{n}^{-k} N_{k}^{\left(n\right)} \right|>{\varepsilon B\left(n\right)\mathord{\left/ {\vphantom {\varepsilon B\left(n\right) 2}} \right. \kern-\nulldelimiterspace} 2} \right\}\left|\mathfrak{F}_{k-1}^{\left(n\right)} \right. \right).$$
Thus, from Chebyshev inequality, for any $\gamma >0$, one has
\begin{equation} \label{eq37}
{\mathbb{P}} \left(L_{1,2} \left(n\right)>\gamma \right)\le \frac{b_{n}}{\gamma B^{2} \left(n\right)} \sum _{k=1}^{\left[nt\right]}a_{n}^{-2k}{\mathbb{E}}X_{k-1}^{\left(n\right)}{\mathbb{E}}\left(I\left\{\left|a_{n}^{-k} N_{k}^{\left(n\right)} \right|>{\varepsilon B\left(n\right)\mathord{\left/ {\vphantom {\varepsilon B\left(n\right) 2}} \right. \kern-\nulldelimiterspace} 2} \right\}\right).
\end{equation}
In order to bound the right-hand side of \eqref{eq37}, we first apply Cauchy-Schwarz inequality, the inequality $\sqrt{\left(x^2+y^2\right)}\leq \sqrt{2}\left(\left|x \right|+\left|y \right| \right)$ and then Chebyshev inequality, overall we obtain
$$ \frac{b_{n} }{\gamma B^{2} \left(n\right)} \sum _{k=1}^{\left[nt\right]}a_{n}^{-2k}\left({\mathbb{E}}\left(X_{k-1}^{\left(n\right)}\right)^{2}\right)^{{1 \mathord{\left/ {\vphantom {1 2}} \right.
\kern-\nulldelimiterspace} 2}} \left({\mathbb{E}}\left[I\left\{\left|a_{n}^{-k} N_{k}^{\left(n\right)} \right|>{\varepsilon B\left(n\right)\mathord{\left/ {\vphantom {\varepsilon B\left(n\right) 2}} \right. \kern-\nulldelimiterspace} 2} \right\}\right]\right)^{{1 \mathord{\left/ {\vphantom {1 2}} \right.
\kern-\nulldelimiterspace} 2}} $$
$$\leq \frac{b_{n}}{\gamma B^{2} \left(n\right)} \sum _{k=1}^{\left[nt\right]}a_{n}^{-2k}\left(B_{n}^{2}\left(k-1\right)+A_{n}^{2}\left(k-1\right)\right)^{{1 \mathord{\left/ {\vphantom {1 2}} \right.
\kern-\nulldelimiterspace} 2}} \left(\mathbb{P}\left\{\left|a_{n}^{-k} N_{k}^{\left(n\right)} \right|>{\varepsilon B\left(n\right)\mathord{\left/ {\vphantom {\varepsilon B\left(n\right) 2}} \right. \kern-\nulldelimiterspace} 2} \right\}\right)^{{1 \mathord{\left/ {\vphantom {1 2}} \right.
\kern-\nulldelimiterspace} 2}}  $$
$$\leq \frac{Cb_{n}}{B^{3} \left(n\right)} \sum_{k=1}^{\left[nt\right]}a_{n}^{-2k}\left(B_{n}\left(k-1\right)+A_{n}\left(k-1\right)\right) \left(\mathbb{E}\left(a_{n}^{-k} N_{k}^{\left(n\right)}\right)^{2} \right)^{{1 \mathord{\left/ {\vphantom {1 2}} \right.
\kern-\nulldelimiterspace} 2}} $$
$$ \leq \frac{Cb_{n}}{B^{3} \left(n\right)} \sum_{k=1}^{\left[nt\right]}a_{n}^{-3k}\beta^{{1 \mathord{\left/ {\vphantom {1 2}} \right.
\kern-\nulldelimiterspace} 2}}\left(n,k\right)\left(B_{n}\left(k-1\right)+A_{n}\left(k-1\right)\right) $$
$$ \leq \frac{Cb_{n}\left(B_{n}\left([nt]\right)+A_{n}\left([nt]\right)\right)}{B^{3} \left(n\right)} \sum_{k=1}^{\left[nt\right]}a_{n}^{-3k}\beta^{{1 \mathord{\left/ {\vphantom {1 2}} \right.
\kern-\nulldelimiterspace} 2}}\left(n,k\right) \to 0,\ \ n\to \infty. $$
where in last step we used the properties of regularly varying functions and $B^2\left(n\right)\sim Cn^{2}\alpha\left(n\right)b_{n}$,
$A\left(n\right)\sim Cn\alpha \left(n\right)$, $n \to \infty$.
Therefore, we deduce that
\begin{equation}\label{eq38}
L_{1,2} \left(n\right)\mathop{\to }\limits^{P} 0, \ \ n\to \infty.
\end{equation}
Recalling \eqref{eq036} and \eqref{eq38}, we have
\begin{equation}\label{eq41}
L_{1} \left(n\right)\mathop{\to }\limits^{P} 0, \ \ n\to \infty.
\end{equation}
In order to estimate $L_{2} \left(n\right)$ it suffices to note that
$$ {\mathbb{E}}L_{2} \left(n\right)\leq \frac{1}{B^{2} \left(n\right)} \sum_{k=1}^{\left[nt\right]}a_{n}^{-2k} {\mathbb{E}} \left(N_{k}^{\left(n\right)} \right)^{2} $$
$$ \leq \frac{2}{B^{2} \left(n\right)} \sum_{k=1}^{\left[nt\right]}a_{n}^{-2k}\left(\beta\left(n,k\right)+{\mathbb{E}} \left({\mathbb{E}} \left(\varepsilon_{k}^{\left(n\right)} \left|\mathfrak{F}_{k-1}^{\left(n\right)} \right. \right)-\alpha\left(n,k\right)\right)^{2}\right) $$
$$ \leq \frac{4}{B^{2} \left(n\right)} \sum_{k=1}^{\left[nt\right]}a_{n}^{-2k}\beta\left(n,k\right) \to 0, \ \ n \to \infty.$$
Therefore,
\begin{equation}\label{eq42}
L_{2} \left(n\right)\mathop{\to }\limits^{P} 0, \ \ n\to \infty.
\end{equation}
From \eqref{eq41}-\eqref{eq42}, one has
\begin{equation}\label{eq43}
L \left(n\right)\mathop{\to }\limits^{P} 0, \ \ n\to \infty.
\end{equation}
Thus, collecting \eqref{eq035} and \eqref{eq43}, we have proved that $Z_{n}^{\left(1\right)} \left(t\right)\mathop{\to }\limits^{D} W\left(\varphi^{\ast}\left(t\right)\right)$, $n\to \infty$ in Skorokhod space $D\left({\mathbb R}_{+},{\mathbb R}\right)$. Now it remains to note that $CW\left(t/C^{2} \right)$ is a standard Wiener process for any $C\neq 0$. This ends the proof of Proposition \ref{Proposition1}.

\begin{proposition}\label{Proposition2}
Assume for each $n\in {\mathbb N}$, $\left\{\varepsilon_{k}^{\left(n\right)} ,\, k\in {\mathbb N}\right\}$ is a sequence of $m-$dependent random variables. If conditions (C1)-(C3) and (C5) hold, then
$$ Z_{n}^{\left(2\right)} \left(t\right)\, \mathop{\to }\limits^{D} \, \, 0,\ \  n\to \infty $$
in Skorokhod space $D\left({\mathbb R}_{+},{\mathbb R}\right)$.\\
\end{proposition}
\textbf{Proof of Proposition \ref{Proposition2}.}
The assertion of Proposition \ref{Proposition2} follows from Lemmas \ref{Lemma2}-\ref{Lemma3} immediately.

Now we are able to prove Theorem \ref{thm2}.\\
\textbf{Proof of Theorem \ref{thm2}.}
The proof of Theorem \ref{thm2} is a straightforward consequence of Propositions \ref{Proposition1}-\ref{Proposition2}.

\section{Appendix}
In the proofs we need the following FLT from Jacod and Shiryaev (\cite{book 3}, Theorem VIII.3.33).
\begin{theorem}\label{thm3}
Let for each $n\ge 1$, $\left\{U_{k}^{n} ,k\ge 1\right\}$ be a sequence of martingale difference with respect to some filtration $\left\{\Re _{k}^{n} ,k\ge 1\right\}$, such that the conditional Lindeberg condition:
\begin{equation} \label{eq39}
\sum_{k=1}^{\left[nt\right]}{\mathbb{E}} \left(\left(U_{k}^{n} \right)^{2} I\left\{\left|U_{k}^{n} \right|>\varepsilon \right\}\left|\Re  _{k-1}^{\left(n\right)} \right. \right) \mathop{\to }\limits^{P} 0, \ \ n\to \infty
\end{equation}
holds for all $\varepsilon >0$ and $t\in {\mathbb R}_{+}$. Then
$$ \sum_{k=1}^{\left[nt\right]}U_{k}^{n}  \mathop{\to }\limits^{D} U\left(t\right),\ \ n\to \infty $$
in Skorokhod space $D\left({\mathbb R}_{+} ,{\mathbb R}\right)$, where $U\left(t\right)$ is a continuous Gaussian martingale with mean zero and covariance function $C\left(t\right)$, $t\in {\mathbb R}_{+}$, if and only if
\begin{equation} \label{eq40}
\sum_{k=1}^{\left[nt\right]}{\mathbb{E}} \left(\left(U_{k}^{n} \right)^{2} \left|\Re_{k-1}^{n} \right. \right) \mathop{\to}\limits^{P} C\left(t\right), \ \  n\to \infty
\end{equation}
for all $t\in {\mathbb R}_{+}$.
\end{theorem}
Now we provide the covariance inequality for random variables satisfying $\psi-$mixing condition (\cite{book 4}, Lemma 1.2.11).
\begin{lemma}\label{Lemma8}
Let $\left\{\xi_{n},n\geq 1\right\}$ be a sequence of  $\psi-$mixing random variables and let $\xi\in \mathfrak{F}_{1}^{k}$,  $\eta \in \mathfrak{F}_{k+n}^{\infty}$, $\mathbb{E}\left|\xi \right|<\infty$, $\mathbb{E}\left|\eta \right|<\infty$. Then
$$ \left|\mathbb{E}\xi\eta - \mathbb{E}\xi \mathbb{E}\eta \right| \leq \psi\left(n\right)\mathbb{E}\left|\xi \right|\mathbb{E}\left|\eta \right|.$$
\end{lemma}

\textbf{Acknowledgements} The author would like to thank professors I.Rahimov and Y.Khusanbaev for many valuable comments which improved the presentation of the paper.


\end{document}